\documentclass{svmult}

\usepackage{makeidx}         
\usepackage{graphicx}        
\usepackage{multicol}        
\usepackage[bottom]{footmisc}
\usepackage{amssymb,amsmath}
\usepackage{longtable}

\def\RR{\mathbb{R}}
\def\SS{\mathbb{S}}
\def\svec{\mbox{svec}}
\def\pennon{{\sc Pennon}{}}

\def\penbmi{{\sc Penbmi}~}
\def\pensdp{{\sc Pensdp}~}
\def\pennon{{\sc Pennon}~}
\def\matlab{{\sc Matlab}~}
\def\yalmip{{\sc Yalmip}~}
\def\ampl{{\sc Ampl}~}

\def\pensdpd{{\sc Pensdp}{}}
\def\pennond{{\sc Pennon}{}}

\def\yalmipd{{\sc Yalmip}{}}

\def\Tr{\mbox{Tr\,}}
\def\COMPlib{{\sl COMP{\sl l$_e\!$ib}} }
\def\real{\mathbb{R}}

\newcommand{\rnn}[2]{\mbox{$\real^{#1 \times #2}$}}

\spdefaulttheorem{Algorithm}{Algorithm}{\bf}{}

\makeindex             

\begin{document}

\title*{PENNON: Software for linear and nonlinear matrix inequalities}
\author{Michal Ko\v{c}vara\inst{1}\and
Michael Stingl\inst{2}}
\institute{School of Mathematics, University of Birmingham, Birmingham
B15 2TT, UK and Institute of Information Theory and Automation, Academy
of Sciences of the Czech Republic, Pod vod\'arenskou
v\v{e}\v{z}\'{\i}~4, 18208 Praha 8, Czech Republic
\texttt{kocvara@maths.bham.ac.uk} \and Institute of Applied
Mathematics, University of Erlangen-Nuremberg, Martensstrasse 3, 91058
Erlangen, Germany \texttt{stingl@am.uni-erlangen.de}}
%
%
\maketitle

\section{Introduction}
The goal of this paper is to present an overview of the software
collection for the solution of linear and nonlinear semidefinite
optimization problems \pennond. In the first part we present theoretical and practical details of
the underlying algorithm and several implementation issues. In the
second part we introduce the particular codes \pensdpd, \penbmi and
\pennond, focus on some specific features of these codes and show
how they can be used for the solution of selected problems.

We use standard notation: $\SS^{m}$ is the space of real symmetric
matrices of dimension $m\times m$ and $\SS^{m}_+$ the space of positive
semidefinite matrices from $\SS^{m}$. The inner product on $\SS^{m}$ is
defined by $\langle A, B\rangle_{\SS^{m}} := {\rm trace} (AB)$.
Notation $A\preccurlyeq B$ for $A,B\in\SS^{m}$ means that the matrix
$B-A$ is positive semidefinite. The norm $\|\cdot\|$ is always the
$\ell_2$ norm in case of vectors and the spectral norm in case of
matrices, unless stated otherwise. Finally, for $\Phi:\SS^m\to\SS^m$
and $X,Y\in \SS^m$, $D\Phi(X)[Y]$ denotes the directional derivative of
$\Phi$ with respect to $X$ in direction $Y$.

\section{The main algorithm}
\subsection{Problem formulation}
The nonlinear semidefinite problems can be written in several different
ways. In this section, for the sake of simplicity, we will use the
following formulation:
\begin{align}
& \min_{x\in\RR^n}  f(x)\label{eq:SDP}\\
  &  \mbox{subject to}\nonumber\\
    &\qquad  {\cal A} (x)  \preccurlyeq  0 \,. \nonumber
\end{align}
Here $f:\RR^n \to \RR$ and ${\cal A}:\RR^n \to \SS^{m}$ are twice
continuously differentiable mappings.

Later in sections on linear SDP, BMI and nonlinear SDP, we will give
more specific formulations of the problem. However, the algorithm and
theory described in this section applies, with some exceptions
discussed later, to all these specific formulations.

\subsection{The algorithm}
The basic algorithm used in this article is based on the nonlinear
rescaling method of R.~Polyak \cite{polyak} and was described in detail
in \cite{pennon} and \cite{stingl}. Here we briefly recall it and
stress points that will be needed in the rest of the paper.

The algorithm is based on the choice of a smooth penalty/barrier
function $\Phi_p:\SS^{m}\to\SS^{m}$ that satisfies a number of
assumptions (see \cite{pennon,stingl}) guaranteeing, in particular,
that for any $p>0$
$$
{\cal A} (x) \preccurlyeq 0 \Longleftrightarrow  \Phi_p({\cal
A}(x)) \preccurlyeq 0
$$
for (at least) all $x$ such that ${\cal A} (x) \preccurlyeq 0$. Thus
for any $p>0$, problem (\ref{eq:SDP}) has the same solution as the
following ``augmented" problem
\begin{align}
  & \min_{x\in\RR^n}  f(x)\label{eq:SDP-phi}\\
  &  \mbox{subject to}\nonumber\\
  &\qquad   \Phi_p({\cal A}(x))  \preccurlyeq  0 \,.\nonumber
\end{align}

The Lagrangian of (\ref{eq:SDP-phi}) can be viewed as a (generalized)
augmented Lagrangian of~(\ref{eq:SDP}):
\begin{equation}\label{eq:lagr}
  F(x,U,p) = f(x)
  + \langle U, \Phi_p \left({\cal
A}(x)\right)\rangle_{\SS^{m}}\,;
\end{equation}
here  $U\in\SS^{m}_+$ is a Lagrangian multiplier associated with the
inequality constraint.

The algorithm below can be seen as a generalization of the Augmented
Lagrangian method.
\begin{Algorithm}\label{algo:1}
Let $x^1$ and $U^1$ be given. Let $p^1>0$.
For $k=1,2,\ldots$ repeat until a stopping criterion is reached:
\begin{enumerate}
\item\quad $x^{k+1}  =  \arg\!\min\limits_{x \in
    \RR^n}F(x,U^k,p^k)$
\item\quad $U^{k+1}  =  D\Phi_p({\cal A}(x^{k+1}))[U^k]$
\item\quad $p^{k+1}  \leq  p^k \,.$
\end{enumerate}
\end{Algorithm}

Details of the algorithm, the choice of the penalty function $\Phi_p$,
the choice of initial values of $x,U$ and $p$, the approximate
minimization in Step~(i) and the update formulas, will be discussed in
subsequent sections. The next section concerns the overview of
theoretical properties of the algorithm.

\subsection{Convergence theory overview}
Throughout this section we make the following assumptions on problem
(\ref{eq:SDP}):
\begin{enumerate}
\item[$(A1)$] $x^*=\arg\min\{f(x)|x\in\Omega\}$ exists, where
    $\Omega=\{x\in\RR^n|{\cal A} (x) \preccurlyeq  0\}$.
\item[$(A2)$] The Karush-Kuhn-Tucker necessary optimality
    conditions hold in $x^*$, i.e., there exists $U^*\in \SS^m$
    such that
\begin{align}
f'(x^*) + \left[\left< U^*,{\cal A}_i\right>\right]_{i=1}^m &= 0 \nonumber \\
\langle U^*, {\cal A}(x^*) \rangle &= 0 \nonumber \\
U^* &\succeq 0  \nonumber \\
{\cal A}(x^*) &\preceq 0,\label{eq:KKT}
\end{align}
where ${\cal A}_i$ denotes the $i$-th partial derivative of ${\cal
A}$ at $x^*$ $(i=1,\ldots,n)$. Moreover the strict complementary is
satisfied.
\item[$(A3)$] The nondegeneracy condition holds, i.e., if for
    $1\leq r < m$ the vectors $s_{m-r+1},\ldots,s_{m} \in \RR^m$
    form a basis of the null space of the matrix ${\cal A}(x^*)$,
    then the following set of $n$-dimensional vectors is linearly
    independent:
$$
v_{i,j}=(s_i^\top {\cal A}_1s_j,\ldots,s_i^\top {\cal
A}_ns_j)^\top,\qquad m-r+1\leq i \leq\ j \leq m.
$$
%
\item[$(A4)$] Define $E_0=(s_{m\!-r\!+1},\ldots,s_m)$, where
    $s_{m\!-r\!+1},\ldots,s_m$ are the vectors introduced in
    assumption $(A3)$. Then the
\emph{cone of critical directions} 
at $x^*$ is defined as
$$
{\cal C}(x^*) =\left\{h\in\RR^n:\sum_{i=1}^n h_i E_0^\top {\cal A}_iE_0
\preceq 0,f'(x^*)^\top h=0\right\}.
$$
With this the following second order sufficient optimality
condition is assumed to hold at $(x^*,U^*)$:~ For all $h\in {\cal
C}(x^*)$ with $h\neq 0$ the inequality
$$
h^\top \left( L''_{xx}(x^*,U^*)+H(x^*,U^*)\right) h > 0,
$$
is satisfied, where $L$ is the classic Lagrangian of (\ref{eq:SDP})
defined as
$$
L(x,U)=f(x) + \left< U,{\cal A}(x)\right>,
$$
$H(x^*,U^*)$ is defined entry-wise by
\begin{equation}
H(x^*,U^*)_{i,j}=-2\left<U^*, {\cal A}_i [{\cal A}(x^*)]^{\dagger}
{\cal A}_j\right> \label{eq:Hstar}
\end{equation} (see, for example, \cite[p. 490]{bonnans-shapiro}) and $[{\cal A}(x^*)]^{\dagger}$ is the
Moore-Penrose inverse of ${\cal A}(x^*)$.
\item[(A5)] Let
$$
\Omega_p=\left\{x\in\RR^n|{\cal A}(x) \preceq pI_{m}\right\}\,.
$$
Then the following growth condition holds:
\begin{equation}
\quad \exists \pi > 0 \mbox{ and } \tau > 0 \mbox{ such that }
\max\left\{\left\|{\cal A}(x)\right\| \mid x\in\Omega_\pi\right\} \leq
\tau\,.  \label{ass:level,noncvx}
\end{equation}
\end{enumerate}

Using these assumptions the following local convergence result can be
established.

\begin{theorem}\label{th:Conv1} Let ${\cal A}(x)$ be twice continuously differentiable and assumptions $(A1)$ to $(A5)$ hold for the pair $(x^*, U^*)$. Then there exists a penalty parameter $p_0>0$ large enough and a neighbourhood ${\cal V}$ of $(x^*,U^*)$ such that for all $(U,p) \in {\cal V}$:
\begin{enumerate}
\item[a)] There exists a vector $$\hat{x} = \hat{x}(U,p) =
    \arg\min\{F(x,U,p)|x\in\RR^n\}$$ such that $\nabla_x
    F(\hat{x},U,p) = 0$.
\item[b)] For the pair $\hat{x}$ and $\widehat{U} =
    \widehat{U}(U,p) = D\Phi_p\left({\cal
    A}(\hat{x}(U,p))\right)\left[U\right]$ the estimate
\begin{equation}
\max\left\{\|\hat{x}-x^*\|, \|\widehat{U}-U^*\| \right\} \leq
Cp\left\|U-U^*\right\|
\end{equation}
holds, where $C$ is a constant independent of $p$.
\item[c)] $\hat{x}(U^*,p)= x^*$ and~ $\widehat{U}(U^*,p)=U^*$.
\item[d)]
The function $F(x,U,p)$ is strongly convex with respect to $x$ in a
neighborhood of $\hat{x}(U,p)$.
\end{enumerate}
\end{theorem}

The proof for Theorem \ref{th:Conv1} as well as a precise definition of
the neighborhood ${\cal V}$ is given in \cite{stingl}. A slightly
modified version of Theorem \ref{th:Conv1} with a particular choice of
the penalty function $\Phi$ can be found in \cite{li-zhang}. An
alternative convergence theorem using slightly different assumptions is
presented in \cite{alpen}.

An immediate consequence of Theorem \ref{th:Conv1} is that Algorithm
\ref{algo:1} converges with a linear rate of convergence. If $p_k\to 0$
for $k\to \infty$ is assumed for the sequence of penalty parameters
then the rate of convergence is superlinear.

\remark a)~ Let $x^+$ be a local minimum of problem (\ref{eq:SDP})
satisfying assumptions (A2) to (A5) and denote by $U^+$ the
corresponding (unique) optimal multiplier. Assume further that there
exists a neighborhood $S_\nu$ of $x^+$ such that there is no further
first order critical point $\tilde{x} \neq x^+$ in $S_\nu(x^+)$. 
Then all statements of Theorem \ref{th:Conv1} remain valid, if we
replace $(x^*, U^*)$ by $(x^+,U^+)$ and the function $\hat{x}(U,p)$ by
\begin{equation} \label{eq:local0}
\hat{x}_{\rm loc}(U,p) = \arg \min\{F(x,U,p)|x\in \RR^n, x \in S_\nu
\}.
\end{equation}
Moreover Theorem \ref{th:Conv1} d) guarantees that $F(x,U,p)$ is
strongly convex in a neighborhood of $x^+$ for all $(U,p) \in {\cal V}
$. Consequently any local descent method applied to the problem
\begin{equation} \label{eq:local}
(i')\qquad \mbox{Find } x^{k+1} \mbox{ such that } \left\|\nabla_x F(x,U^k,p^k)\right\| = 0\\
\end{equation}
will automatically find a solution, which satisfies the additional
constraint $x^{k+1}\in S_\nu$ provided it is started with $x^k$
close enough to $x^+$. Moreover, Algorithm 1 will converge to the local optimum $x^+$ (see \cite{stingl} for more details). %

\medskip\noindent
b)~ A global convergence result can be found in \cite{stingl}.

\subsection{{Choice of $\Phi_p$}}\label{sec:hess} The penalty function
$\Phi_p$ of our choice is defined as follows:
\begin{equation}\label{eq:pen}
\Phi_p({\cal A}(x))  = -p^2({\cal A}(x) - pI)^{-1} - pI \,.
\end{equation}
The advantage of this choice is that it gives closed formulas for the
first and second derivatives of $\Phi_p$. Defining
\begin{equation}\label{eq:Z}
  {\cal Z}(x) = -({\cal A}(x) - pI)^{-1}
\end{equation}
we have (see \cite{pennon}):
\begin{align}
\frac{\partial}{\partial x_i} \Phi_p({\cal A}(x))& =
p^2{\cal Z}(x) \frac{\partial{\cal A}(x)}{\partial x_i} {\cal Z}(x)
\label{eq:der1} \\
\frac{\partial^2}{\partial x_i\partial x_j} \Phi_p({\cal
A}(x)) & = p^2{\cal Z}(x) \left(\frac{\partial{\cal A}(x)}{\partial
x_i}
{\cal Z}(x) \frac{\partial{\cal A}(x)}{\partial x_j} +
 \frac{\partial^2{\cal A}(x)}{\partial x_i\partial x_j}
  \right.\nonumber\\
&\left.\phantom{p^2{\cal Z}(x)}\qquad\ \
+ \frac{\partial{\cal A}(x)}{\partial x_j}
{\cal Z}(x) \frac{\partial{\cal A}(x)}{\partial x_i}\right){\cal
Z}(x)\,.&
\label{eq:der2}
\end{align}

\subsection{The modified Newton method} To solve the (possibly
nonconvex) unconstrained minimization problem in Step 1, we use the
following modification of the Newton method with line-search:
\begin{Algorithm} \label{algo:5}
Given an initial iterate $x_0$, repeat for all $k=0,1,2,\ldots$ until a
stopping criterion is reached
\begin{enumerate}
\item Compute the gradient~$g_k$ and Hessian~$H_k$ of $F$ at~$x_k$.
\item Try to factorize $H_k$ by Cholesky decomposition. If $H_k$ is
    factorizable, set $\widehat{H} = H_k$ and go to Step~4. \item
    Compute $\beta \in \left[-\lambda_{\rm min}, -2\lambda_{\rm
    min}\right]$, where $\lambda_{\rm min}$ is the minimal
    eigenvalue of~$H_k$ and set
$ \widehat{H} = H_k + \beta I. $
\item Compute the search direction $ d_k = -\widehat{H}^{-1}g_k. $
\item Perform line-search in direction~$d_k$. Denote the
    step-length by~$s_k$.
\item Set $ x_{k+1} =  x_k + s_k d_k. $
\end{enumerate}
\end{Algorithm}
The step-length $s$ in direction $d$ is calculated by a gradient free
line-search algorithm that tries to satisfy the Armijo condition.
Obviously, for a convex $F$, Algorithm \ref{algo:5} is just the damped
Newton method, which is known to converge under standard assumptions.

If, in the non-convex case, the Cholesky factorization in Step~2 fails,
we calculate the value of~$\beta$ in Step~3 in the following way:
\begin{Algorithm} \label{algo:6}
For a given $\beta_0 > 0$
\begin{enumerate}
\item Set $\beta = \beta_0$.
\item Try to factorize $H + \beta I$ by the Cholesky method.
\item If the factorization fails due to a negative pivot element,
    go to step 4, otherwise go to step 5.
\item If $\beta \geq \beta_0$, set $\beta= 2{\beta}$ and continue
    with 2. Otherwise go to step 6.
\item If $\beta \leq \beta_0$, set $\beta= \frac{\beta}{2}$ and continue with step 2. Otherwise STOP.
\item Set $\beta = 2\beta$ and STOP.
\end{enumerate}
\end{Algorithm}
Obviously, when Algorithm \ref{algo:6} terminates we have $\beta \in
\left[-\lambda_{\rm min}, -2\lambda_{\rm min}\right]$. It is well known
from the nonlinear programming literature that under quite mild assumptions any cluster point
of the sequence generated by Algorithm \ref{algo:5} is a first order
critical point of problem in Step~1 of Algorithm~1.
\begin{remark}
There is one exception, when we use a different strategy for the
calculation of $\beta$. The exception is motivated by the observation
that the quality of the search direction gets poor, if we choose
$\beta$ too close to $-\lambda_{\rm min}$. Therefore,
 if we encounter bad quality of the search direction, we
use a bisection technique to calculate an approximation of
$\lambda_{\rm min}$, denoted by $\lambda_{\rm min}^{\rm a}$, and
replace $\beta$ by $-1.5\lambda_{\rm min}^{\rm a}$.
\end{remark}
\begin{remark}
Whenever we will speak about the Newton method or Newton system, later
in the paper, we will always have in mind the modified method described
above.
\end{remark}

\subsection{{How to solve the linear systems?}}\label{sec:LinSys}
In both algorithms proposed in the preceding sections one has to solve
repeatedly linear systems of the form
\begin{equation} \label{eq:linsys2}
(H+D) d = -g,
\end{equation}
where $D$ is a diagonal matrix chosen such that the matrix $H+D$ is
positive definite. There are two categories of methods, which can be
used to solve problems of type (\ref{eq:linsys2}): direct and iterative
methods. Let us first concentrate on the direct methods.

\subsubsection{Cholesky method}
Since the system matrix in (\ref{eq:linsys2}) is always positive
definite, our method of choice is the Cholesky method. Depending on the
sparsity structure of $H$, we use two different realizations:
\begin{itemize}
\item If the fill-in of the Hessian is below $20 \%$\,, we use a
    sparse Cholesky solver which is based on ideas of Ng and Peyton
    \cite{ng-peyton}. The solver makes use of the fact that the
    sparsity structure is the same in each Newton step in all
    iterations. Hence the sparsity pattern of $H$, reordering of
    rows and columns to reduce the fill-in in the Cholesky factor,
    and symbolic factorization of $H$ are all performed just once
    at the beginning of Algorithm 1.
Then, each time the system (\ref{eq:linsys2}) has to be solved, the
numeric factorization is calculated based on the precalculated
symbolic factorization. Note that we added stabilization techniques
described in \cite{wright} to make the solver more robust for
almost singular system matrices.
\item Otherwise, if the Hessian is dense, we use the {\sc atlas}
    implementation of the {\sc lapack} Cholesky solver {\tt
    DPOTRF}.
\end{itemize}

\subsubsection{Iterative methods} We solve the system $\widehat{H}d = -g$ with a symmetric positive
definite and, possibly, ill-conditioned matrix $\widehat{H}=H+D$. We
use the very standard preconditioned conjugate gradient method. The
algorithm is well known and we will not repeat it here. The algorithm
is stopped when the normalized residuum is sufficiently small:
$$
  \|\widehat{H}d_k+g\|/\|g\| \leq \epsilon \,.
$$
In our tests, the choice $\epsilon = 5\cdot 10^{-2}$ was sufficient.

\subsubsection{Preconditioners} We are looking for a preconditioner---a
matrix $M\in\SS^n_+$---such that the system $M^{-1}\widehat{H}d =
-M^{-1}g$ can be solved more efficiently than the original system
$\widehat{H}d=-g$. Apart from standard requirements that the
preconditioner should be efficient and inexpensive, we also require
that it should only use Hessian-vector products. This is particularly
important in the case when we want to use the Hessian-free version of
the algorithm.

\paragraph{Diagonal preconditioner}
This is a simple and often-used preconditioner with
$$
  M={\rm diag}\,(\widehat{H}).
$$
On the other hand, being simple and general, it is not considered to be
very efficient. Furthermore,  we need to know the diagonal elements of
the
Hessian. It is certainly possible to compute 
these elements 
by Hessian-vector products. For that, however, we would need $n$
gradient evaluations and the approach would become too costly.

\paragraph{L-BFGS preconditioner}
Introduced by Morales-Nocedal \cite{morales-nocedal}, this
preconditioner is intended for application within the Newton method.
(In a slightly different context, the \mbox{(L-)BFGS} preconditioner
was also proposed in \cite{fukuda}.) The algorithm is based on the
limited-memory BFGS formula (\cite{nocedal-wright}) applied to
successive CG (instead of Newton) iterations. The preconditioner, as
used in \pennon is described in detail in \cite{pen-iter}. Here we only
point out some important features.

As recommended in the standard L-BFGS method, we used 16--32 correction
pairs, if they were available. Often the CG method finished in less
iterations and in that case we could only use the available iterations
for the correction pairs. If the number of CG iterations is higher than
the required number of correction pairs $\mu$, we may ask how to select
these pairs. We have two options: Either we take the last $\mu$ pairs
or an ``equidistant'' distribution over all CG iterations. The second
option is slightly more complicated but we may expect it to deliver
better results.

The L-BFGS preconditioner has the big advantage that it only needs
Hessian-vector products and can thus be used in the Hessian-free
approaches. On the other hand, it is more complex than the above
preconditioners; also our results are not conclusive concerning the
efficiency of this approach. For many problems it worked
satisfactorily, for some, on the other hand, it even lead to higher
number of CG steps than without preconditioner.

\subsection{Multiplier and penalty update}\label{sec:mult}
For the penalty function $\Phi_p$ from (\ref{eq:pen}), the formula for
update of the matrix multiplier $U$ in Step (ii) of
Algorithm~\ref{algo:1} reduces to
\begin{equation}\label{eq:upd}
U^{k+1} = (p^k)^2{\cal Z}(x^{k+1}) U^k {\cal Z}(x^{k+1})
\end{equation}
with ${\cal Z}$ defined as in (\ref{eq:Z}). Note that when $U^k$ is
positive definite, so is $U^{k+1}$. We set $U^1$ equal to a positive
multiple of the identity.

Numerical tests indicate that big changes in the multipliers should be
avoided for the following reasons. Big change of $U$ means big change
of the augmented Lagrangian that may lead to a large number of Newton
steps in the subsequent iteration. It may also happen that already
after few initial steps the multipliers become ill-conditioned and the
algorithm suffers from numerical difficulties. To overcome these, we do
the following:
\begin{enumerate}
\item Calculate $U^{k+1}$ using (\ref{eq:upd}).
\item Choose a positive $\mu_{\! A} \leq 1$, typically 0.5.
\item Compute $\lambda_{\! A}=\min\left(\mu_{\! A}, \mu_{\! A}
    \frac {\left\|U^k\right\|_F}
    {\left\|U^{k+1}-U^k\right\|_F}\right).$
\item Update the current multiplier by
$$
U^{new} = U^k + \lambda_{\! A} (U^{k+1} - U^k).
$$
\end{enumerate}

Given an initial iterate $x^1$, the initial penalty parameter $p^1$ is
chosen large enough to satisfy the inequality
$$
p^1 I - {\cal A}(x^1) \succ 0.
$$
Let $\lambda_{\rm max}({\cal A}(x^{k+1})) \in \left(-\infty,p^k\right)$
denote the maximal eigenvalue of ${\cal A}(x^{k+1})$, $\pi < 1$ be a
constant factor, depending on the initial penalty parameter $p^1$
(typically chosen between $0.3$ and $0.6$) and $x_{\rm{feas}}$ be a
feasible point. Let $l$ be set to 0 at the beginning of Algorithm \ref
{algo:1}. Using these quantities, our strategy for the penalty
parameter update can be described as follows:

\begin{enumerate}
\item If $p^k < p_{eps}$, set $\gamma = 1$ and go to 6.
\item Calculate $\lambda_{\rm max}({\cal A}(x^{k+1}))$.
\item If $\pi p^{k} > \lambda_{\rm max}({\cal A}(x^{k+1}))$, set
    $\gamma = \pi,$ $l=0$ and go to 6.
\item If $l < 3$, set $\gamma = \left(\lambda_{\rm max}({\cal
    A}(x^{k+1})) + p^k\right)/(2p^k)$, set $l:=l+1$ and go to 6.
\item Let $\gamma = \pi$, find $\lambda \in \left(0,1\right)$ such,
    that $$\lambda_{\rm max} \left({\cal A}(\lambda x^{k+1} +
    (1-\lambda) x_{\rm{feas}})\right) < \pi p^k,$$ set $x^{k+1} =
    \lambda x^{k+1} + (1-\lambda) x_{\rm{feas}}$ and $l:=0$.
\item Update current penalty parameter by $ p^{k+1} = \gamma p^k. $
\end{enumerate}
The reasoning behind steps 3 to 5 is as follows: As long as the
inequality
\begin{equation} \label{eq:penup}
\lambda_{\rm max}({\cal A}(x^{k+1})) < \pi p^k
\end{equation}
holds, the values of the augmented Lagrangian in the next iteration
remain finite and we can reduce the penalty parameter by the predefined
factor $\pi$. As soon as inequality (\ref{eq:penup}) is violated, an
update using $\pi$ would result in an infinite value of the augmented
Lagrangian in the next iteration. Therefore the new penalty parameter
should be chosen from the interval $(\lambda_{\rm max}({\cal
A}(x^{k+1})), p^k)$. Because a choice close to the left boundary of the
interval leads to large values of the augmented Lagrangian, while a
choice close to the right boundary slows down the algorithm, we choose
$\gamma$ such that
$$
p^{k+1} = \frac{\lambda_{\rm max}({\cal A}(x^{k+1})) + p^k}{2} \,.
$$
In order to avoid stagnation of the penalty parameter update process
due to repeated evaluations of step 4, we redefine $x^{k+1}$ using the
feasible point $x_{\rm{feas}}$ whenever step 4 is executed in three
successive iterations; this is controlled by the parameter $l$. If no
feasible point is yet available, Algorithm \ref {algo:1} is stopped and
restarted from the scratch with a different choice of initial
multipliers. The parameter $p_{eps}$ is typically chosen as $10^{-6}$.
In case we detect problems with convergence of Algorithm \ref {algo:1},
$p_{eps}$ is decreased and the penalty parameter is updated again,
until the new lower bound is reached.

\subsection{Initialization and stopping criteria}\label{sec:stopping}
\paragraph{Initialization}
Algorithm 1 can start with an arbitrary primal variable $x \in \RR^n$.
Therefore we simply choose $x^1 = 0$. For the description of the
multiplier initialization strategy we rewrite problem (SDP) in the
following form:
\begin{align}
& \min_{x\in\RR^n}  f(x)\label{eq:bSDP}\\
&
  \mbox{subject to}\nonumber \\ &\qquad
   {\cal A}_i (x)  \preccurlyeq  0, \quad \,i=1,\ldots,\sigma\,.\nonumber
\end{align}
Here ${\cal A}_i (x)\in\SS^{m_j}$ are diagonal blocks of the original
constrained matrix ${\cal A} (x)$ and we have $\sigma=1$ if ${\cal A}
(x)$ consists of only one block. Now the initial values of the
multipliers are set to
$$
\begin{aligned}
U^1_j &= \mu_j I_{m_j}, &\quad j=1,\ldots,\sigma,
\end{aligned}
$$
where $I_{m_j}$ are identity matrices of order $m_j$ and
\begin{align}
\mu_j &= m_j\max_{1\leq \ell\leq n}\frac{1+\left|\frac{\partial f(x)}{\partial x_\ell}\right|}
  {1+\left\|\frac{\partial{\cal A}(x)}{\partial x_\ell}\right\|}\,.
\end{align}
Given the initial iterate $x^1$, the initial penalty parameter $p^1$ is
chosen large enough to satisfy the inequality
$$
p^1 I - {\cal A}(x^1) \succ 0.
$$
\paragraph{Stopping criterion in the sub-problem}  In the first
iterations of Algorithm \ref{algo:1}, the approximate minimization of
$F$ is stopped when $\|\frac{\partial}{\partial x} F(x,U,p)\| \leq
\alpha$, where $\alpha = 0.01$ is a good choice in most cases. In the
remaining iterations, after a certain precision is reached, $\alpha$ is
reduced in each outer iteration by a constant factor, until a certain
$\underline{\alpha}$ (typically $10^{-7}$) is reached.
\paragraph{Stopping criterion for the main algorithm}
We have implemented two different stopping criteria for the main
algorithm.
\begin{itemize}
\item \emph{First alternative:} The main algorithm is stopped if
    both of the following inequalities hold:
\begin{eqnarray*}
\frac{|f(x^k) - F(x^k,U^k,p)|}{1+|f(x^k)|} < \varepsilon_1 \,, \quad
\quad \frac{|f(x^k) - f(x^{k-1})|}{1+|f(x^k)|} < \varepsilon_1 \,,
\end{eqnarray*}
where $\varepsilon_1$ is typically $10^{-7}$.
\item \emph{Second alternative:} The second stopping criterion is
    based on the KKT-conditions. Here the algorithm is stopped,
    if$$ \min\left\{\lambda_{\max} \left({\cal
A}(x)\right), \left|\langle {\cal A}(x), U\rangle\right|,
\|\nabla_x F(x,U,p)\| \right\} \leq \varepsilon_2.
$$
\end{itemize}

\subsection{Complexity}
The computational complexity of Algorithm 1 is clearly dominated by
Step 1. In each step of the Newton method, there are two critical
issues: assembling of the Hessian of the augmented Lagrangian and
solution of the linear system of equations (the Newton system). 

\subsubsection{Hessian assembling}
\paragraph{Full matrices}
Assume first that all the data matrices are full. The assembling of the
Hessian (\ref{eq:der2}) can be divided into the following steps:
\begin{itemize}
\item Calculation of ${\cal Z}(x)  \longrightarrow \, O(m^3+m^2n)$.
\item Calculation of ${\cal Z}(x) U {\cal Z}(x) \longrightarrow \,
    O(m^3)$.
\item Calculation of ${\cal Z}(x) U {\cal Z}(x){\cal A}'_i(x) {\cal
    Z}(x)$ for all $i$ $\, \longrightarrow \,O(m^3n)$.
\item Assembling the rest $\, \longrightarrow O(m^2 n^2)$.
\end{itemize}
Now it is straightforward to see that an estimate of the complexity of
assembling of (\ref{eq:der2}) is given by $O(m^3n + m^2n^2)$.

Many optimization problems, however, have very sparse data structure
and therefore have to be treated by sparse linear algebra routines. We
distinguish three basic types of sparsity.

\paragraph{The block diagonal case}
The first case under consideration is the block diagonal case. In
particular, we want to describe the case, where\\[-0.5cm]
\begin{itemize}
\item the matrix
${\cal A}(x)$ consists of many (small) blocks.\\[-0.5cm]
\end{itemize}
In this situation the original SDP problem (\ref{eq:SDP}) can be
written in the form (\ref{eq:bSDP}). If we define
$\bar{m}=\max\{m_i~|~i=1,\ldots,d\}$ we can estimate the computational
complexity of the Hessian assembling by $O(d\bar{m}^3n+d\bar{m}^2n^2)$.
An
interesting subcase of problem (\ref{eq:bSDP}) is when\\[-0.5cm]
\begin{itemize}
\item each of the matrix constraints ${\cal A}_i (x)$ involves
just a few components of $x$.\\[-0.5cm]
\end{itemize}
If we denote the maximal number of components of $x$ on which each of
the blocks ${\cal A}_i(x), i=1,2,\ldots,d$ depends by $\bar{n}$, our
complexity formula becomes $O(d\bar{m}^3\bar{n}+d\bar{m}^2\bar{n}^2)$.
If we further assume that the numbers $\bar{n}$ and $\bar{m}$ are of
order $O(1)$, then the complexity estimate can be further simplified to
$O(d)$. A typical example for this sparsity class are the `mater'
problems discussed in section \ref{sec:pensdp}.

\paragraph{The case when ${\cal A}(x)$ is dense and ${\cal A}'_i(x)$ are sparse}
Let us first mention that for any index pair $(i,j) \in
\{1,\ldots,n\}\times\{1,\ldots,n\}$ the non-zero structure of the
matrix ${\cal A}''_{i,j}(x)$ is given by (a subset of the) intersection
of the non-zero index sets of the matrices ${\cal
A}'_i(x)$ and ${\cal A}'_j(x)$. We assume that\\[-0.5cm]
\begin{itemize}
\item there are at most $O(1)$ non-zero entries in ${\cal A}'_i(x)$
    for
all $i=1,\ldots,n$.\\[-0.5cm]
\end{itemize}
Then the calculation of the term $$ \left[\left\langle {\cal Z}(x)
  U {\cal Z}(x),
     {\cal A}''_{i,j}(x) \right\rangle\right]_{i,j = 1}^n$$
can be performed in $O(n^2)$ time. In the paper by Fujisawa, Kojima and
Nakata on exploiting sparsity in semidefinite programming
\cite{fujisawa-kojima-nakata} several ways are presented how to
calculate a matrix of the form
\begin{equation} \label{eq:matmult}
D_1 S_1 D_2 S_2
\end{equation}
efficiently, if $D_1$ and $D_2$ are dense and $S_1$ and $S_2$ are
sparse matrices. If our assumption above holds, the calculation of the
matrix
$$
\left[\left\langle {\cal Z}(x) U {\cal Z}(x)
     {\cal A}'_j(x)
     {\cal Z}(x), {\cal A}'_i(x)
     \right\rangle\right]_{i,j = 1}^n
$$
can be performed in $O(n^2)$ time. Thus, recalling that for the
calculation of ${\cal Z}(x)$ we have to compute the inverse of an
$(m\!\times \! m)$-matrix, we get the following complexity estimate for
the Hessian assembling: $O(m^3+n^2)$. Note that in our implementation
we follow the ideas presented in \cite{fujisawa-kojima-nakata}. Many
linear SDP problems coming from real world applications have exactly
the sparsity structure discussed in this paragraph.

\paragraph{The case when ${\cal A}(x)$ and the Cholesky factor of ${\cal A}(x)$ is sparse}
Also in this case we can conclude that all partial derivatives of
${\cal A}(x)$ of first and second order are sparse matrices. Therefore
it suffices to assume that
\\[-0.5cm]
\begin{itemize}
\item the matrix ${\cal A}(x)$ has at most $O(1)$ non-zero entries.\\[-0.5cm]
\end{itemize}
We have to compute expressions of type
$$
({\cal A}(x) - pI)^{-1} U ({\cal A}(x) - pI)^{-1} \quad \mbox{ and }
\quad ({\cal A}(x) - pI)^{-1}.
$$
Note that each of the matrices above can be calculated by maximally two
operations of the type $(A - I)^{-1}M$, where $M$ is a symmetric
matrix. Now assume that not only ${\cal A}(x)$ but also its Cholesky
factor is sparse. Then, obviously, the Cholesky factor of $({\cal
A}(x)-pI)$, denoted by $L$, will also be sparse.
This leads to the following assumption:\\[-0.5cm]
\begin{itemize}
\item Each column of $L$ has at most $O(1)$ non-zero entries.\\[-0.5cm]
\end{itemize}
Now the $i$-th column of $C:=({\cal A}(x) - pI)^{-1} M$ can then be
computed as
$$
  C^i = (L^{-1})^TL^{-1}M^i,\ i=1,\ldots,n,
$$
and the complexity of computing $C$ by Cholesky factorization is
$O(m^2)$, compared to $O(m^3)$ when computing the inverse of $(A(x) -
pI)$ and its multiplication by $U$. So the overall complexity of
Hessian assembling  is of order $O(m^2+n^2)$.

\begin{remark}
Recall that in certain cases, we do not need to assemble the Hessian
matrix. In this case the complexity estimates can be improved
significantly; see the next section.
\end{remark}

\subsubsection{Solution of the Newton system}
As mentioned above, the Newton system
\begin{equation}
  Hd = -g \label{eq:system}
\end{equation}
can either be solved by a direct (Cholesky) solver or by an iterative
method.

\paragraph{Cholesky method}
The complexity of Cholesky algorithm is $ O(n^3) $ for dense matrices
and $O(n^{\kappa})$, $1\leq\kappa\leq 3$ for sparse matrices, where
$\kappa$ depends on the sparsity structure of the matrix, going from a
diagonal to a full matrix.

\paragraph{Iterative algorithms}
From the complexity viewpoint, the only demanding step in the CG method
is a matrix-vector product with a matrix of dimension $n$. For a dense
matrix and vector, it needs $O(n^2)$ operations. Theoretically, in
exact arithmetics, the CG method needs $n$ iterations to find an exact
solution of the system, hence it is equally expensive as the Cholesky
algorithm. There are, however, two points that may favor the CG method.

First, it is well known that the convergence behavior of the CG method
can be significantly improved by preconditioning. The choice of the
preconditioner $M$ will be the subject of the next section.

The second---and very important---point is that we actually do not need
an exact solution of the Newton system. On the contrary, a rough
approximation of it will do (see \cite[Thm. 10.2]{kanzow}). Hence, in
practice, we may need just a few CG iterations to reach the required
accuracy. This is in contrast with the Cholesky method where we cannot
control the accuracy of the solution and always have to compute the
exact one (within the machine precision). Note that we always start the
CG method with initial approximation $d_0=0$; thus, performing just one
CG step, we would obtain the steepest descend method. Doing more steps,
we improve the search direction toward the Newton direction; note the
similarity to the Toint-Steihaug method \cite{nocedal-wright}.

Summarizing these two points: when using the CG algorithm, we may
expect to need just $O(n^2)$ operations, at least for well-conditioned
(or well-preconditioned) systems.

Note that we are still talking about dense problems. The use of the CG
method is a bit nonstandard in this context---usually it is preferable
for large sparse problems. However, due to the fact that we just need a
very rough approximation of the solution, we may favor it to the
Cholesky method also for medium-sized dense problems.

\paragraph{Approximate Hessian formula}
When solving the Newton system by the CG method, the Hessian is only
needed in a matrix-vector product of the type $Hv := \nabla^2 F(x^k)v
$. Because we only need to compute the products, we may use a finite
difference formula for the approximation of this product
\begin{equation}\label{eq:approx}
  \nabla^2 F(x^k)v \approx \frac{\nabla F(x^k+hv) - \nabla F(x^k)}{h}
\end{equation}
with $h = (1+\|x^k\|_2\sqrt\varepsilon)$; see \cite{nocedal-wright}. In
general, $\varepsilon$ is chosen so that the formula is as accurate as
possible and still not influenced by round-off errors. The ``best''
choice is obviously case dependent; in our implementation, we use
$\varepsilon=10^{-6}$. Hence the complexity of the CG method amounts to
the number of CG iterations times the complexity of gradient
evaluation, which is of order $O(m^3 + Kn)$, where $K$ denotes the
maximal number of nonzero entries in $\mathcal{A}'_i(x),
i=1,2,\ldots,n$. This may be in sharp contrast with the Cholesky method
approach when we have to compute the full Hessian \emph{and} solve the
system by Cholesky method. Again, we have the advantage that we do not
have to store the Hessian in the memory.

This approach is clearly not always applicable. With certain SDP
problems it may happen that the Hessian computation is not much more
expensive than the gradient evaluation. In this case the Hessian-free
approach may be rather time-consuming. Indeed, when the problem is
ill-conditioned and we need many CG iterations, we have to evaluate the
gradient many (thousand) times. On the other hand, when using Cholesky
method, we compute the Hessian just once.


\section{PENSDP}\label{sec:pensdp}
When both functions in (\ref{eq:SDP}) are linear, the problem
simplifies to a standard (primal or dual, as you like) linear
semidefinite programming problem (LSDP)
\begin{align}
 &\min_{x\in\RR^n}  f^T x  \label{eq:LSDP} \\
  &\mbox{subject to} \nonumber\\
 &\qquad  \sum_{k=1}^{n} x_k A_k - A_0 \preccurlyeq 0\,. \nonumber
\end{align}
Here we write explicitly all matrix inequality constraints, as well as
linear constraints, in order to introduce necessary notation. In the
next sections, we will present some special features of the code
\pensdp designed to solve (\ref{eq:LSDP}), as well as selected
numerical examples demonstrating it capabilities.

\subsection{The code PENSDP}
\subsubsection{Special features}
\paragraph{Stopping criteria}
In the case of linear semidefinite programs, we have additionally
adopted the DIMACS criteria \cite{dimacs1}. To define these criteria,
we denote $\widetilde{\cal A} (x) = \sum_{k=1}^{n} x_k A_k$. Recall
that $U$ is the corresponding Lagrangian multiplier and let
$\widetilde{\cal A}^*(\cdot)$ denote the adjoint operator to
$\widetilde{\cal A}(\cdot)$. The DIMACS error measures are defined as
\begin{alignat*}{2}
 \mbox{err}_1 &= \frac{\|\widetilde{\cal A}^*(U) - f\|}{1+\|f\|}\\
  \mbox{err}_2 &= \max\left\{0,\frac{-\lambda_{\rm min}
  (U)}{1+\|f\|}\right\} &\qquad
  \mbox{err}_4 &= \max\left\{0,\frac{-\lambda_{\rm min} (\widetilde{\cal A}(x)-A_0)}{1+\|A_0\|}\right\}  \\
  \mbox{err}_5 &= \frac{\langle A_0,U\rangle - f^Tx}{1+|\langle
  A_0,U\rangle| +|f^Tx|} &\qquad \mbox{err}_6 &= \frac{\langle
  \widetilde{\cal A}(x)-A_0,U\rangle }{1+|\langle A_0,U\rangle| +|f^Tx|}  \,.
\end{alignat*}
Here, $\mbox{err}_1$ represents the (scaled) norm of the gradient of
the Lagrangian, $\mbox{err}_2$ and $\mbox{err}_4$ is the dual and
primal infeasibility, respectively, and $\mbox{err}_5$ and
$\mbox{err}_6$ measure the duality gap and the complementarity
slackness. Note that, in our code, $\mbox{err}_2=0$ by definition; also
$\mbox{err}_3$ that involves the slack variable (not used in our
problem formulation) is automatically zero. If the ``DIMACS stopping
criterion'' is activated we require that
$$
  \mbox{err}_k \leq \delta_{\scriptscriptstyle\rm DIMACS}, \quad
  k\in\{1,4,5,6\} \,.
$$
\paragraph{Implicit Hessian formula}
As mentioned before, when solving the Newton system by the CG method,
the Hessian is only needed in a matrix-vector product of the type $Hv
:= \nabla^2 F(x^k)v $. Instead of computing the Hessian matrix
explicitly and then multiplying it by a vector $v$, we can use the
following formula for the Hessian-vector multiplication
\begin{equation}\label{eq:implicit}
\nabla^2 F(x^k) v = 2{\cal A}^* \left((p^k)^2 {\cal Z}(x^k) U^k {\cal
Z}(x^k) {\cal A}(v) {\cal Z}(x^k)\right),
\end{equation}
where we assume that ${\cal A}$ is linear of the form $A(x) =
\sum_{i=1}^n x_i A_i$ and ${\cal A}^*$ denotes its adjoint. Hence, in
each CG step, we only have to evaluate matrices ${\cal A}(v)$ (which is
simple), ${\cal Z}(x^k)$ and ${\cal Z}(x^k) U^k {\cal Z}(x^k)$ (which
are needed in the gradient computation, anyway), and perform two
additional matrix-matrix products. The resulting complexity formula for
one Hessian-vector product is thus $O(m^3 + Kn)$, where again $K$
denotes the maximal number of nonzero entries in $\mathcal{A}'_i(x),
i=1,2,\ldots,n$.

The additional (perhaps the main) advantage of this approach is the
fact that we do not have to store the Hessian in the memory, thus the
memory requirements (often the real bottleneck of SDP codes) are
drastically reduced.
\paragraph{Dense versus sparse}
For the efficiency of \pensdpd, it is important to know if the problem
has a sparse or dense Hessian. The program can check this
automatically. The check, however, may take some time and memory, so if
the user knows that the Hessian is dense (and this is the case of most
problems), this check can be avoided. This, for certain problems, can
lead to substantial savings not only in CPU time but also in memory
requirements.

\paragraph{Hybrid mode}
For linear semidefinite programming problems, we use the following
hybrid approach, whenever the number of variables $n$ is large compared
to the size of the matrix constraint $m$: We try to solve the linear
systems using the iterative approach as long as the iterative solver
needs a moderate number of iterations. In our current implementation
the maximal number of CG iterations allowed is 100. Each time the
maximal number of steps is reached, we solve the system again by the
Cholesky method. Once the system is solved by the Cholesky method, we
use the Cholesky factor as a preconditioner for the iterative solver in
the next system. As soon as the iterative solver fails three times in
sequel, we completely switch to the Cholesky method.

The hybrid mode allows us to reach a high precision solution while
keeping the solution time low. The main reason is that, when using the
iterative approach, the Hessian of the Augmented Lagrangian has not to
be calculated explicitly.

\subsubsection{User interfaces}
The user has a choice of several interfaces to \pensdpd.
\paragraph{{\sc sdpa} interface}
The problem data are written in an ASCII input file in a {\sc sdpa}
sparse format, as introduced in \cite{sdpa}. The code needs an
additional ASCII input file with parameter values.

\paragraph{{\sc C}/{\sc C++}/{\sc Fortran} interface}%
\pensdp can also be called as a function (or subroutine) from a {\sc
C}, {\sc C++} or {\sc Fortran} program. In this case, the user should
link the \pensdp library to his/her program. In the program the user
then has to specify problem dimensions, code parameters and the problem
data (vectors and matrices) in a sparse format.

\paragraph{{\sc Matlab} interface}
In {\sc Matlab}, \pensdp is called with the following arguments:

\begin{verbatim}
[f,x,u,iflag,niter,feas] = pensdpm(pen);
\end{verbatim}
where {\tt pen} a {\sc Matlab} structure array with fields describing
the problem dimensions and problem data, again in a sparse format.

\paragraph{{\sc Yalmip} interface}
The most comfortable way of preparing the data and calling \pensdp is
via \yalmip \cite{yalmip}. \yalmip is a modelling language for advanced
modeling and solution of convex and nonconvex optimization problems. It
is implemented as a free toolbox for MATLAB. When calling \pensdp from
\yalmip, the user does not has to  bother with the sparsity pattern of
the problem---any linear optimization problem with vector or matrix
variables will be translated by \yalmip into \pensdp data structure.

\subsection{Numerical experiments}
It is not our goal to compare \pensdp with other linear SDP solvers.
This is done elsewhere in this book and the reader can also consult the
benchmark page of Hans
Mittelmann\footnote{plato.la.asu.edu/bench.html}, containing
contemporary results. We will thus present only results for selected
problems and will concentrate on the effect of special features
available in \pensdpd. The results for the `mater' and `rose13'
problems were obtained on an Intel Core i7 processor 2.67GHz with 4GB
memory.

\subsubsection{Sparsity: `mater' problems}
Let us consider the `mater*' problems from Mittelmann's
collection\footnote{plato.asu.edu/ftp/sparse\_sdp.html}. These problems
are significant by several different sparsity patterns of the problem
data. The problem has many small matrix constraints, the data matrices
are sparse, only very few variables are involved in each constraint and
the resulting Hessian matrix is sparse. We cannot switch off sparsity
handling in routines for Hessian assembling but we can run the code
with (forced use of) dense Cholesky factorization and with sparse
Cholesky routine. For instance, problem `mater3' with 1439 variables
and 328 matrix constraints of size 11 was solved in 32 seconds using
the dense Cholesky and only 4 seconds using the sparse Cholesky
routine. The difference is, of course, more dramatic for larger
problems. The next problem `mater4' has 4807 variables and 1138 matrix
constraints of size 11. While the sparse version of \pensdp only needed
20 seconds to solve it, the dense version needed 1149 second. And while
the largest problem `mater6' (20463 variables and 4968 matrix
constraints) does not even fit in the 4GB memory for the dense version,
the sparse code needs only 100MB and solves the problem in 134 seconds.

\subsubsection{Iterative solver: `TOH' collection}
The effect of the use of preconditioned conjugate gradient method for
the solution of the Newton system was described in detail in
\cite{pen-iter,pen-itere}. Recall that iterative solvers are suitable
for problems with a large number of variable and relatively small
constraint matrices. We select from \cite{pen-iter,pen-itere} two
examples arising from maximum clique problems on randomly generated
graphs (the 'TOH' collection in \cite{pen-iter}). The first example is
`theta62' with 13390 variables and matrix size 300. This problem could
still be solved using the direct (dense Cholesky) solver and the code
needed 13714 seconds to solve it. Compared to that, the iterative
version of the code only needed 40 seconds to obtain the solution with
the same precision. The average number of CG steps in each Newton
system was only 10. The largest problem solved in the paper was
`theta162' with 127600 variables and a matrix constraint of size 800.
Note that the Hessians of this example is \emph{dense}, so to solve the
problem by the direct version of \pensdp (or by any other
interior-point algorithm) would mean to store and factorize a full
matrix of dimension 127600 by 127600. On the other hand, the iterative
version of \pensdpd, being effectively a first-order code, has only
modest memory requirements and allowed us to solve this problem in only
672 seconds.

\subsubsection{Hybrid mode: `rose13'}
To illustrate the advantages of the hybrid mode, we consider the
problem `rose13' from Mittelmann's
collection\footnote{plato.asu.edu/ftp/sparse\_sdp.html}. The problem
has 2379 variables and one matrix constraint of size 105. When we solve
the problem by \pensdp with a direct solver of the Newton system, the
code needs 17 global iterations, 112 Newton steps and the solution is
obtained in 188 seconds CPU time, 152 seconds of which is spent in the
Cholesky factorization routine.

Let us now solve the problem using the iterative solver for the Newton
systems. Below we see the first and the last iterations of \pensdpd.
The required precision of DIMACS criteria is
$\delta_{\scriptscriptstyle\rm DIMACS} = 10^{-3}$.
\begin{verbatim}
****************************************************
* it |     obj      |     opt      |   Nwt |  CG   *
****************************************************
|   0|  0.0000e+000 |  0.0000e+000 |     0 |     0 |
|   1|  1.8893e+003 |  8.3896e+000 |    10 |   321 |
|   2|  2.2529e+002 |  8.2785e+000 |    17 |  1244 |
...
|   9| -1.1941e+001 |  2.2966e+000 |    36 |  9712 |
|  10| -1.1952e+001 |  4.9578e+000 |    46 | 10209 |
...
|  15| -1.1999e+001 |  5.0429e-002 |   119 | 103905 |
|  16| -1.1999e+001 |  4.4050e-003 |   134 | 167186 |
****************************************************
\end{verbatim}
The table shows the global iterations of Algorithm 1, the value of the
objective function and the gradient of the augmented Lagrangian and, in
the last two columns, the cumulative number of Newton and CG steps. The
code needed a large number of CG steps that was growing with increasing
conditioning of the Newton system. The problem was solved in 732
seconds of CPU time. When we try to solve the same problem with a
higher precision ($\delta_{\scriptscriptstyle\rm DIMACS} =10^{-7}$),
the iterative method, and consequently the whole algorithm, will get
into increasing difficulties. Below we see the last two iterations of
\pensdp before it was stopped due to one-hour time limit.
\begin{verbatim}
...
|  28| -1.2000e+001 |  5.2644e-002 |   373 | 700549 |
|  29| -1.2000e+001 |  3.0833e-003 |   398 | 811921 |
****************************************************
\end{verbatim}
We can see that the optimality criterium is actually oscillating around
$10^{-3}$.

We now switch the hybrid mode on. The difference will be seen already
in the early iterations of \pensdpd. Running the problem with
$\delta_{\scriptscriptstyle\rm DIMACS}=10^{-3}$, we get the following
output
\begin{verbatim}
****************************************************
* it |     obj      |     opt      |   Nwt |  CG   *
****************************************************
|   0|  0.0000e+000 |  0.0000e+000 |      0 |    0 |
|   1|  1.8893e+003 |  8.3896e+000 |     10 |  321 |
|   2|  2.3285e+002 |  1.9814e+000 |     18 |  848 |
...
|   9| -1.1971e+001 |  4.5469e-001 |     36 | 1660 |
|  10| -1.1931e+001 |  7.4920e-002 |     63 | 2940 |
...
|  13| -1.1998e+001 |  4.1400e-005 |    104 | 5073 |
|  14| -1.1999e+001 |  5.9165e-004 |    115 | 5518 |
****************************************************
\end{verbatim}
The CPU time needed was only 157 seconds, 130 of which were spent in
the Cholesky factorization routine. When we now increase the precision
to $\delta_{\scriptscriptstyle\rm DIMACS}=10^{-7}$, the \pensdp with
the hybrid mode will only need a few more iterations to reach it:
\begin{verbatim}
...
|  16| -1.2000e+001 |  9.8736e-009 |    142 | 6623 |
|  17| -1.2000e+001 |  5.9130e-007 |    156 | 7294 |
****************************************************
\end{verbatim}
The total CPU time increased to 201 seconds, 176 of which were spent in
Cholesky factorization.

Notice that the difference between the direct solver and the hybrid
method would be even more significant for larger problems, such as
`rose15'.


\section{PENBMI}\label{sec:penbmi}
We solve the SDP problem with quadratic objective function and linear
and bilinear matrix inequality constraints:
\begin{align}\label{eq:bmi}
 &\min_{x\in\RR^n} \frac{1}{2} x^T Q x + f^T x \\
 & \mbox{subject to}\nonumber\\
  &\qquad \sum_{k=1}^n b^i_k x_k \leq c^i, \quad i=1,\ldots,N_\ell \nonumber\\
 &\qquad A^i_0 + \sum_{k=1}^{n} x_k A^i_k + \sum_{k=1}^{n}\sum_{\ell=1}^{n} x_k x_\ell K^i_{k\ell}
 \preccurlyeq 0, \quad i=1,\ldots,N\,, \nonumber
\end{align}
where all data matrices are from $\SS^m$.

\subsection{The code PENBMI}
%
\subsubsection{User interface}
The advantage of our formulation of the BMI problem is that, although
nonlinear, the data only consist of matrices and vectors, just like in
the linear SDP case. The user does not have to provide the (first and
second) derivatives of the matrix functions, as these are readily
available. Hence the user interface of \penbmi is a direct extension of
the interface to \pensdp described in the previous section.

In particular, the user has the choice of calling \penbmi from Matlab
or from a C/C++/Fortran code. In both cases, the user has to specify
the matrices $A^i_k, k=1,\ldots n,$ and $K^i_{k\ell}, k,\ell=1,\ldots
n,$ for all constraints $i=1,\ldots,N$, matrix $Q$ from the objective
function and vectors $f, c, b^i, i=1,\ldots,N$. As in the linear SDP
case, all matrices and vectors are assumed to be sparse (or even void),
so the user has to provide the sparsity pattern of the constraints
(which matrices are present) and sparsity structure of each matrix.

Again, the most comfortable way of preparing the data and calling
\penbmi is via \yalmipd. In this case, the user does not has to stick
to the formulation (\ref{eq:bmi}) and bother with the sparsity pattern
of the problem---any optimization problem with vector or matrix
variables and linear or quadratic objective function and (matrix)
constraints will be translated by \yalmip into formulation
(\ref{eq:bmi}) and the corresponding user interface will be
automatically created. Below is a simple example of \yalmip code for
the LQ optimal feedback problem formulated as
\begin{align*}
 &\min_{P\in\RR^{2\times 2},\,K\in\RR^{1\times 2}} \mbox{trace} (P) \\
 & \begin{aligned}
  \mbox{s.t.}\qquad\quad
  &(A+BK)^TP+P(A+BK) \prec -I_{2\times 2}-K^TK\\
 &P\succ 0
 \end{aligned}
\end{align*}
with
$$
  A = \begin{pmatrix}-1& 2\\-3& -4\end{pmatrix},\qquad B = \begin{pmatrix}1\\1\end{pmatrix}\,.
$$
Using \yalmipd, the problem is formulated and solved by the following
few lines
\begin{verbatim}
>> A = [-1 2;-3 -4]; B = [1;1];
>> P = sdpvar(2,2); K = sdpvar(1,2);
>> F = [P >= 0; (A+B*K)'*P+P*(A+B*K) <= -eye(2)-K'*K];
>> optimize(F,trace(P),sdpsettings('solver','penbmi'));
\end{verbatim}

\subsection{{The Static Output Feedback Problem}} \label{sec:static}
Many interesting problems in linear and nonlinear systems control
cannot be formulated and solved as LSDP. BMI formulation of the control
problems was made popular in the mid 1990s \cite{goh94}; there were,
however, no computational methods for solving non-convex BMIs, in
contrast with convex LMIs for which powerful interior-point algorithms
were available.

The most fundamental of these problems is perhaps static output
feedback (SOF) stabilization: given a triplet of matrices $A$,$B$,$C$
of suitable dimensions, find a matrix $F$ such that the eigenvalues of
matrix $A+BFC$ are all in a given region of the complex plane, say the
open left half-plane \cite{blondel00}.

No LSDP formulation is known for this problem but a straightforward
application of Lyapunov's stability theory leads to a BMI formulation:
matrix $A+BFC$ has all its eigenvalues in the open left half-plane if
and only if there exists a matrix $X$ such that
$$
(A+BFC)^T X+(A+BFC)X \prec 0, \quad X=X^T \succ 0
$$
where  $\prec 0$ and $\succ 0$ stand for positive and negative
definite, respectively.

We present a short description of the benchmark collection
\mbox{\COMPlib:} the {\sl CO}nstrained {\sl M}atrix--optimization {\sl
P}roblem {\sl lib}rary
\cite{leibfritz:COMPleib:COnstrained:Matrix:Opt:Prob:lib}\footnote{See
{\tt http://www.mathematik.uni-trier.de/$\sim$leibfritz/}\\{\tt
Proj\_TestSet/NSDPTestSet.htm}}. \COMPlib can be used as a benchmark
collection for a very wide variety of algorithms solving matrix
optimization problems. Currently \COMPlib consists of 124 examples
collected from the engineering literature and real-life applications
for LTI control systems of the form
\begin{equation} \label{COMPlibStateSpaceControlSystem}
 \begin{array}{rcl}
   \dot{x}(t) &=& A x(t) + B_1 w(t) +  B u(t),  \\
   z(t) &=& C_1 x(t) + D_{11} w(t) + D_{12} u(t), \\
   y(t) &=& C x(t) + D_{21} w(t),
 \end{array}
\end{equation}
where $x \in \real^{n_x}$, $u \in \real^{n_u}$, $y \in \real^{n_y}$, $z
\in \real^{n_z}$, $w \in \real^{n_w}$ denote the state, control input,
measured output, regulated output, and noise input, respectively.

The heart of \COMPlib is the MATLAB function file {\sl COMPleib.m}.
This function returns the data matrices $A$, $B_1$, $B$, $C_1$, $C$,
$D_{11}$, $D_{12}$ and $D_{21}$ of
(\ref{COMPlibStateSpaceControlSystem}) of each individual \COMPlib
example.
Depending on specific control design goals, it is possible to derive
particular matrix optimization problems using the data matrices
provided by {\sl COMP{\sl l$_e\!$ib}}. A non exhaustive list of matrix
optimization problems arising in feedback control design are stated in
\cite{leibfritz:COMPleib:COnstrained:Matrix:Opt:Prob:lib}. Many more
control problems leading to NSDPs, BMIs or SDPs can be found in the
literature.

Here we state the BMI formulation of two basic static output feedback
control design problems: SOF--${\cal H}_2$ and SOF--${\cal
H}_{\infty}$. The goal is to determine the matrix $F \in
\rnn{n_u}{n_y}$ of the SOF control law $u(t)=Fy(t)$ such that the
closed loop system
\begin{equation} \label{ClosedLoopSystem}
  \begin{array}{rcl}
   \dot{x}(t) &=& A(F) x(t) +  B(F) w(t),  \\
   z(t) &=& C(F)x(t)+D(F) w(t),
 \end{array} 
\end{equation}
fulfills some specific control design requirements, where $A(F)=A+BFC$,
$B(F)=B_1+BFD_{21}$, $C(F)=C_1+D_{12}FC$, $D(F)=D_{11}+D_{12}FD_{21}$.

We begin with the SOF--${\cal H}_2$ problem: {\em Suppose that
$D_{11}=0$ and $D_{21}=0$. Find a SOF gain $F$ such that $A(F)$ is
Hurwitz and the ${\cal H}_2$--norm of (\ref{ClosedLoopSystem}) is
minimal.} This problem can be rewritten to the following ${\cal
H}_2$--BMI problem formulation, see,
e.g.~\cite{leibfritz:COMPleib:COnstrained:Matrix:Opt:Prob:lib}:
\begin{equation} \label{OptimalH2BMI}
   \begin{array}{c}
     \min \;\; Tr(X)\quad
     \mbox{s.t.~} \;\; Q \succ 0, \\[0.2cm]
       (A+BFC)Q+Q(A+BFC)^T+B_1B_1^T \preceq 0,
     \\[2ex]
       \left [ \begin{array}{cc} X  & (C_1+D_{12}FC)Q \\ Q(C_1+D_{12}FC)^T & Q \end{array}
        \right ] \succeq 0,
   \end{array}
\end{equation}
where $Q \in \rnn{n_x}{n_x}$, $X \in \rnn{n_z}{n_z}$.

${\cal H}_{\infty}$ synthesis is an attractive model--based control
design tool and it allows incorporation of model uncertainties in the
control design. The optimal SOF--${\cal H}_{\infty}$ problem can be
formally stated in the following term: {\em Find a SOF matrix $F$ such
that $A(F)$ is Hurwitz and
 the ${\cal H}_{\infty}$--norm of (\ref{ClosedLoopSystem}) is
minimal.} We consider the following well known ${\cal H}_{\infty}$--BMI
version, see,
e.g.~\cite{leibfritz:COMPleib:COnstrained:Matrix:Opt:Prob:lib}:
\begin{equation} \label{AlternativeOptimalHinfBMIformulation}
  \begin{array}{c}
   \min \;\; \gamma \quad
   \mbox{s.t.~} \;\; X \succ 0, \quad \gamma > 0,\\
                  \left [ \begin{array}{ccc}
                    A(F)^TX+XA(F) & XB(F) & C(F)^T \\ B(F)^TX & -\gamma \; I_{n_w} & D(F)^T
                    \\ C(F) & D(F) & -\gamma \; I_{n_z}
                     \end{array} \right ] \prec 0,
  \end{array}
\end{equation}
where $\gamma \in \real$, $X \in \rnn{n_x}{n_x}$.

We present results of our numerical experiences for the static output
feedback problems of \COMPlib. The link between \COMPlib and PENBMI was
provided by the MATLAB parser YALMIP~3 \cite{yalmip}.
All tests were performed on a 2.5\,GHz Pentium with 1\,GB RDRAM under
Linux. The results of PENBMI for ${\cal H}_2$-BMI and ${\cal
H}_{\infty}$-BMI problems can be divided into seven groups: The first
group consists of examples solved without any difficulties (38 problems
in the ${\cal H}_2$ case and 37 problems for the ${\cal H}_{\infty}$
setting). The second and third group contain all cases for which we had
to relax our stopping criterion. In 4 (11) examples the achieved
precision was still close to our predefined stopping criterion, while
in 5 (7) cases deviation is significant (referring to ${\cal H}_2$
(${\cal H}_{\infty}$)). Then there are examples, for which we could
calculate almost feasible solutions, but which failed to satisfy the
Hurwitz-criterion, namely AC5 and NN10. The fourth and fifth group
consist of medium and small scale cases for which {\sc PENBMI} failed,
due to ill conditioned Hessian of $F$---the Cholesky algorithm used for
its factorization did not deliver accurate solution and the Newton
method failed. In the ${\cal H}_2$-setting (fourth group) these are
AC7, AC9, AC13, AC18, JE1, JE2, JE3, REA4, DIS5, WEC1, WEC2, WEC3, UWV,
PAS, NN1, NN3, NN5, NN6, NN7, NN9, NN12 and NN17, in the ${\cal
H}_{\infty}$-setting (fifth group) JE1, JE2, JE3, REA4, DIS5, UWV, PAS,
TF3, NN1, NN3, NN5, NN6, NN7 and NN13. The cases in the sixth group are
large scale, ill conditioned problems, where {\sc PENBMI} ran out of
time (AC10, AC14, CSE2, EB5). Finally, for very large test cases our
code runs out of memory (HS1, BDT2, EB6, TL, CDP, NN18).

\subsection{{Simultaneous Stabilization BMIs}} \label{sec:simstab}

Another example leading to BMI formulation is the problem of
simultaneously stabilizing a family of single-input single-output
linear systems by one fixed controller of given order. This problem
arises for instance when trying to preserve stability of a control
system under the failure of sensors, actuators, or processors.
Simultaneous stabilization of three or more systems was extensively
studied in \cite{blondel94}. Later on, the problem was shown to belong
to the wide range of robust control problems that are
NP-hard~\cite{blondel00}.

In \cite{scl} a BMI formulation of the simultaneous stabilization
problem was obtained in the framework of the polynomial, or algebraic
approach to systems control. This formulation leads to a feasibility
BMI problem which, in a more general setting can be reformulated by the
following procedure: Assume we want to find a feasible point of the
following system of BMIs
\begin{equation} \label{eq:BMIfeas}
  A^i_0 + \sum_{k=1}^{n} x_k A^i_k + \sum_{k=1}^{n}\sum_{\ell=1}^{n} x_k x_\ell K^i_{k\ell} \prec 0, \qquad i=1,\ldots,N
\end{equation}
with symmetric matrices $A_k^i,K^i_{k\ell} \in \RR^{d_i\times d_i}$,
$k,\ell=1,\ldots,n$,  $i=1,\ldots,N$, and $x\in\RR^n$. Then we can
check the feasibility of (\ref{eq:BMIfeas}) by solving the following
optimization problem
\begin{align}\label{eq:BMI}
& \min_{x\in\RR^n,\lambda\in\RR}  \lambda \\
& 
  \mbox{s.t.}\qquad
  A^i_0 + \sum_{k=1}^{n} x_k A^i_k + \sum_{k=1}^{n}\sum_{\ell=1}^{n} x_k x_\ell K^i_{k\ell} \preccurlyeq  \lambda I_n,& \qquad i=1,\ldots,N
  \,.
\end{align}
Problem (\ref{eq:BMI}) is a global optimization problem: we know that
if its global minimum $\lambda$ is non-negative then the original
problem (\ref{eq:BMIfeas}) is infeasible. On the other hand \penbmi can
only find critical points, so when solving~(\ref{eq:BMI}), the only
conclusion we can make is the following:
\begin{verse}
when $\lambda < 0$, the system is strictly feasible;\\
when $\lambda = 0$, the system is marginally feasible;\\
when $\lambda > 0$ the system may be infeasible.
\end{verse}

During numerical experiments it turned out that the feasible region
of~(\ref{eq:BMIfeas}) is often unbounded. We used two strategies to
avoid numerical difficulties in this case: First we introduced large
enough artificial bounds $x_{\rm bound}$. Second, we modify the
objective function by adding the square of the 2-norm of the vector $x$
multiplied by a weighting parameter $w$. After these modifications
problem (\ref {eq:BMI}) reads as follows:
\begin{align}\label{eq:BMIa}
& \min_{x\in\RR^n,\lambda\in\RR}  \lambda + w \|x\|^2_2\\
& \begin{aligned}
  \mbox{s.t.}\qquad\qquad\qquad\qquad
    -x_{\rm bound} \leq x^k \leq x_{\rm bound},& \qquad k=1,\ldots,n \\
  A^i_0 + \sum_{k=1}^{n} x_k A^i_k + \sum_{k=1}^{n}\sum_{\ell=1}^{n} x_k x_\ell K^i_{k\ell} \preccurlyeq  \lambda I_{n\times n},& \qquad i=1,\ldots,N
  \,.
\end{aligned} \notag
\end{align}
This is exactly the problem formulation we used in our numerical
experiments.

Results of numerical examples for a suite of simultaneous stabilization
problems selected from the recent literature can be found in
\cite{penbmi0}.


\section{PENNON}\label{sec:pennon}
\subsection{The problem and the modified algorithm}
\subsubsection{Problem formulation}
In this, so far most general version of the code, we solve optimization
problems with a nonlinear objective subject to nonlinear inequality and
equality constraints and semidefinite bound constraints:
\begin{align}
& \min_{x\in\RR^n, Y_1\in\SS^{p_1},\ldots,Y_k\in\SS^{p_k}}  f(x,Y)\label{eq:nlpsdp}\\
& \begin{aligned}
  \mbox{subject to}\quad
   &g_i(x,Y)  \leq 0, \qquad &&i=1,\ldots,m_g\nonumber\\
   &h_i(x,Y)  = 0, \qquad &&i=1,\ldots,m_h \nonumber\\
   &\underline{\lambda}_i I \preceq Y_i\preceq \overline{\lambda}_i I,
     \qquad&&i=1,\ldots,k \,.\nonumber
\end{aligned}
\end{align}
Here
\begin{itemize}
\item $x\in\RR^n$ is the vector variable
\item $Y_1\in\SS^{p_1},\ldots,Y_k\in\SS^{p_k}$ are the matrix
    variables; we denote $Y=(Y_1,\ldots,Y_k)$
\item $f$, $g_i$ and $h_i$ are $C^2$ functions from $\RR^n\times
    \SS^{p_1}\times\ldots\times\SS^{p_k}$ to $\RR$
\item $\underline{\lambda}_i$ and $\overline{\lambda}_i$ are the
    lower and upper bounds, respectively, on the eigenvalues of
    $Y_i$, $i=1,\ldots,k$
\end{itemize}

Although the semidefinite inequality constraints are of a simple type,
most nonlinear SDP problems can be formulated in the above form. For
instance, the problem (\ref{eq:SDP}) can be transformed into
(\ref{eq:nlpsdp}) using slack variables and equality constraints, when
$$
  {\cal A}(x)\preccurlyeq 0
$$
is replaced by
\begin{align*}
  &{\cal A}(x) = S \quad \mbox{element-wise}\\
  &S\preccurlyeq 0
\end{align*}
with a new matrix variable $S\in\SS^m$.

\subsubsection{Direct equality handling}
Problem (\ref{eq:nlpsdp}) is not actually a problem of type
(\ref{eq:SDP}) that was introduced in the first section and for which
we have developed the convergence theory. The new element here are the
equality constraints. Of course, we can formulate the equalities as two
inequalities, and this works surprisingly well for many problems.
However, to treat the equalities in a ``proper'' way, we adopted a
concept which is successfully used in modern primal-dual interior point
algorithms (see, e.g., \cite{ipopt}): rather than using augmented
Lagrangians, we handle the equality constraints directly on the level
of the subproblem. This leads to the following approach. Consider the
optimization problem
\begin{equation}\label{eq:SDP-eq}
\begin{aligned}
  \qquad \qquad & \min_{x\in\RR^n}  f(x)\\
  &   \mbox{subject to}\\
  &\qquad   {\cal A}(x)  \preccurlyeq  0 \,,\\
  &\qquad   h (x)  =  0 \,,
\end{aligned}
\end{equation}
where $f$ and ${\cal A}$ are defined as in the previous sections and
$h:\RR^n \to \RR^d$ represents a set of equality constraints. Then we
define the augmented Lagrangian
\begin{eqnarray}
  &&\overline{F} (x,U,v,p) =\nonumber \\
  &&f(x)
  + \langle U, \Phi_{p} ( {\cal
A}(x))\rangle_{\SS^{m}}  + v^\top
h(x) \,, \label{eq:lagr3}
\end{eqnarray}
where $U,\Phi,p$ are defined as before and $v\in \RR^d$ is the vector
of Lagrangian multipliers associated with the equality constraints.
Now, on the level of the subproblem, we attempt to find an approximate
solution of the following system (in $x$ and $v$):
\begin{equation}\label{eq:KKT-eq}
\begin{aligned}
  \qquad \nabla_x \overline{F} (x,U,v,p) &= 0\,, \\
 h(x) &= 0\,,
\end{aligned}
\end{equation}
where the penalty parameter $p$ as well as the multiplier $U$ are
fixed. In order to solve systems of type (\ref{eq:KKT-eq}), we apply
the damped Newton method. Descent directions are calculated utilizing
the factorization routine MA27 from the Harwell subroutine library
(\cite{dure:82}) in combination with an inertia correction strategy as
described in \cite{ipopt}. Moreover, the step length is derived using
an augmented Lagrangian merit function defined as
$$
\overline{F} (x,U,v,p) + \frac{1}{2\mu}\|h(x)\|_2^2
$$
along with an Armijo rule.

\subsubsection{Strictly feasible constraints}
In certain applications, the bound constraints must remain strictly
feasible for all iterations because, for instance, the objective
function may be undefined at infeasible points \cite{fmo-stress}. To be
able to solve such problems, we treat these inequalities by a classic
barrier function. For this reason we introduce an additional matrix
inequality
$${\cal
S}(x)\preccurlyeq 0 $$ in problem (\ref{eq:SDP-phi})  and define the
augmented Lagrangian
\begin{equation}\label{eq:lagr2}
  \widetilde{F} (x,U,p,s) = f(x)
  + \langle U, \Phi_{p} ( {\cal
A}(x))\rangle_{\SS^{m}} + s \Phi_{\rm bar}({\cal S}(x)) \,,
\end{equation}
where $\Phi_{\rm bar}$ can be defined, for example, by
$$\Phi_{\rm bar}({\cal
S}(x)) = -\log\det(-{\cal S}(x)).$$
Note that, while the penalty parameter $p$ maybe constant from a
certain index $\bar{k}$ (see again \cite{stingl} for details), the
barrier parameter $s$ is required to tend to zero with increasing $k$.

\subsection{The code PENNON}
\subsubsection{Slack removal}
As already mentioned, to transform constraints of the type ${\cal
A}(x)\preccurlyeq 0 $  into our standard structure, we need to
introduce a slack matrix variable $S$, and replace the original
constraint by ${\cal A}(x) = S$ element-wise, and $S\preccurlyeq 0$.
Thus in order to formulate the problem in the required form, we have to
introduce a new (possibly large) matrix variable and many new equality
constraints, which may have a negative effect on the performance of the
algorithm. However, the reformulation using slack variables is only
needed for the input of the problem, not for its solution by Algorithm
1. Hence, the user has the option to say that certain matrix variables
are actually slacks and these are then automatically removed by a
preprocessor. The code then solves the problem with the original
constraint ${\cal A}(x)\preccurlyeq 0$.

\subsubsection{User interface}
Unlike in the \pensdp and \penbmi case, the user has to provide not
only function values but also the first and second derivatives of the
objective and constraint functions. In the \matlab and C/C++/Fortran
interface the user is required to provide six functions/subroutines for
evaluation of function value, gradient and Hessian of the objective
function and the constraints, respectively, at a given point.

To make things simple, the matrix variables are treated as vectors in
these functions, using the operator $\svec: \SS^m\to\RR^{(m+1)m/2}$
defined by
$$
  \svec\begin{pmatrix}a_{11}&a_{12}&\ldots&a_{1m}\\
                            &a_{22}&\ldots&a_{2m}\\
                            &      &\ddots&\vdots\\
                            sym            &&&a_{mm}
  \end{pmatrix}
  =(a_{11},a_{12},a_{22},\ldots, a_{1m},
                             a_{2m},\ldots,a_{mm})^T
$$
In the main program, the user defines problem sizes, values of bounds,
and information about matrix variables (number, sizes and sparsity
patterns).

In addition, we also provide an interface to \ampl \cite{ampl} which is
a comfortable modelling language for optimization problems. As \ampl
does not support matrix variables, we treat them, within an \ampl
script, as vectors, using the operator $\svec$ defined above.

\begin{example} Assume that we have a matrix variable $X\in\SS^3$
$$
  X = \begin{pmatrix} x_1&x_2&x_4\\ x_2&x_3&x_5\\x_4&x_5&x_6
  \end{pmatrix}
$$
and a constraint
$$
  \Tr (X A) = 3 \qquad\mbox{with}\ A = \begin{pmatrix} 0&0&1\\ 0&1&0\\1&0&0
  \end{pmatrix}\,.
$$
The matrix variable is treated as a vector
$$\svec(X) = (x_1,x_2,\ldots,x_6)^T$$
and the above constraint is thus equivalent to
$$
  x_3+2x_4 = 3\,.
$$
\end{example}

The code needs to identify the matrix variables, their number and size.
These data are included in an ASCII file that is directly read by the
\pennon and that, in addition, includes information about lower and
upper bounds on these variables.

\subsection{Examples}
Most examples of nonlinear semidefinite programs that can be found in
the literature are of the form: for a given (symmetric, indefinite)
matrix $H$ find a nearest positive semidefinite matrix satisfying
possibly some additional constraints. Many of these problems can be
written as follows
\begin{align}
  &\min_{X\in\SS^n} \frac{1}{2}\|X-H\|_F^2\label{eq:sdp_quad}\\
  &\mbox{subject to}\nonumber\\
  &\qquad \langle A_i,X\rangle=b_i,\quad i=1,\ldots,m\nonumber\\
  &\qquad X\succeq 0\nonumber
\end{align}
with $A_i\in\SS^n$, $i=1,\ldots,m$. Probably the most prominent example
is the problem of finding the nearest correlation matrix \cite{higham}.

Several algorithms have been derived for the solution of this problems;
see, e.g., \cite{higham,malick}. It is not our primal goal to compete
with these specialized algorithms (although \pennon can solve problems
of type (\ref{eq:sdp_quad}) rather efficiently). Rather we want to
utilize the full potential of our code and solve ``truly nonlinear''
semidefinite problems. In the rest of this section we will give
examples of such problems.

\subsection{Correlation matrix with the constrained condition number}\label{ex:cond}
We consider the problem of finding the nearest correlation matrix:
\begin{align}
&\min_X \sum_{i,j=1}^n (X_{ij}-H_{ij})^2\label{corr1}\\
&\mbox{subject to}\nonumber\\
&\qquad X_{ii} = 1,\quad i=1,\ldots,n\nonumber\\
&\qquad X\succeq 0\nonumber
\end{align}

We will consider an example based on a practical application from
finances; see \cite{werner-schoettle}. Assume that a $5\times 5$
correlation matrix is extended by one row and column. The new data is
based on a different frequency than the original part of the matrix,
which means that the new matrix is no longer positive definite:
$$
  H_{\rm ext} = \begin{pmatrix}1 &-0.44& -0.20 &0.81& -0.46& -0.05\\
    -0.44& 1 &0.87& -0.38& 0.81 & -0.58\\
    -0.20 &.87 &1& -0.17& 0.65& -0.56\\
    0.81 &-0.38& -0.17& 1& -0.37& -0.15\\
    -0.46& 0.81& 0.65& -0.37& 1& -0.08\\
    -0.05&-0.58&-0.56&-0.15&0.08&1
    \end{pmatrix}\,.
$$
Let us find the nearest correlation matrix to $H_{\rm ext}$ by solving
(\ref{corr1}) (either by \pennon or by any of the specialized
algorithms mentioned at the beginning of this section). We obtain the
following result (for the presentation of results, we will use {\sc
Matlab} output in short precision):
\begin{verbatim}
X =
    1.0000   -0.4420   -0.2000    0.8096   -0.4585   -0.0513
   -0.4420    1.0000    0.8704   -0.3714    0.7798   -0.5549
   -0.2000    0.8704    1.0000   -0.1699    0.6497   -0.5597
    0.8096   -0.3714   -0.1699    1.0000   -0.3766   -0.1445
   -0.4585    0.7798    0.6497   -0.3766    1.0000    0.0608
   -0.0513   -0.5549   -0.5597   -0.1445    0.0608    1.0000
\end{verbatim}
with eigenvalues
\begin{verbatim}
eigen =
    0.0000    0.1163    0.2120    0.7827    1.7132    3.1757
\end{verbatim}
As we can see, one eigenvalue of the nearest correlation matrix is
zero. This is highly undesirable from the application point of view. To
avoid this, we can add lower (and upper) bounds on the matrix variable,
i.e., constraints $
  \underline{\lambda} I \preceq X\preceq \overline{\lambda} I $.
However, the application requires a different approach when we need to
\emph{bound the condition number of the nearest correlation matrix},
i.e., to add the constraint
$$
  \mbox{cond}(X) \leq \kappa \,.
$$
This constraint can be introduced in several ways. For instance, we can
introduce the constraint
$$
  I\preceq\widetilde{X}\preceq \kappa I
$$
using the transformation $
  \widetilde{X} = \zeta X
$. The problem of finding the nearest correlation matrix with a
\emph{given} condition number then reads as follows:
\begin{align}
&\min_{\zeta,\widetilde{X}}
\sum_{i,j=1}^n (\frac{1}{\zeta}\widetilde{X}_{ij}-H_{ij})^2\label{corr_cond}\\
&\mbox{subject to}\nonumber\\
&\qquad  \widetilde{X}_{ii} -\zeta = 0,\quad i=1,\ldots,n\nonumber\\
&\qquad I\preceq\widetilde{X}\preceq \kappa I\,.\nonumber
\end{align}
The new problem now has the NLP-SDP structure of (\ref{eq:nlpsdp}).
When solving it by \pennon with $\kappa=10$, we get the solution after
11 outer and 37 inner iterations. The optimal value of $\zeta$ is
$3.4886$ and, after the back substitution $X =
\frac{1}{\zeta}\widetilde{X}$, we get the nearest correlation matrix
\begin{verbatim}
X =
    1.0000   -0.3775   -0.2230    0.7098   -0.4272   -0.0704
   -0.3775    1.0000    0.6930   -0.3155    0.5998   -0.4218
   -0.2230    0.6930    1.0000   -0.1546    0.5523   -0.4914
    0.7098   -0.3155   -0.1546    1.0000   -0.3857   -0.1294
   -0.4272    0.5998    0.5523   -0.3857    1.0000   -0.0576
   -0.0704   -0.4218   -0.4914   -0.1294   -0.0576    1.0000
\end{verbatim}
with eigenvalues
\begin{verbatim}
eigenvals =
    0.2866    0.2866    0.2867    0.6717    1.6019    2.8664
\end{verbatim}
and the condition number equal to 10, indeed.

\paragraph{Large-scale problems}
To test the capability of the code to solve large-scale problems, we
have generated randomly perturbed correlation matrices $H$ of arbitrary
dimension by the commands
\begin{verbatim}
n = 500; x=10.^[-4:4/(n-1):0];
G = gallery('randcorr',n*x/sum(x));
E = 2*rand(n,n)-ones(n,n); E=triu(E)+triu(E,1)'; E=(E+E')/2;
H = (1-0.1).*G + 0.1*E;
\end{verbatim}
For large-scale problems, we successfully use the iterative
(preconditioned conjugate gradient) solver for the Newton system in
Step 1 of the algorithm. In every Newton step, the iterative solver
needs just a few iterations, making it a very efficient alternative to
a direct solver.

For instance, to solve a problem with a $500\times 500$ matrix $H$ we
needed 11 outer and 148 inner iterations, 962 CG steps, and 21 minutes
on a notebook. Note that the problem had 125251 variables and 500
linear constraints: that means that at each Newton step we solved
(approximately) a system with a \emph{full} $125251\times 125251$
matrix. The iterative solver (needing just matrix-vector product) was
clearly the only alternative here.

We have also successfully solved a real-world problem with a matrix of
dimension 2000 and with many additional linear constraints in about 10
hours on a standard Linux workstation with 4 Intel Core 2 Quad
processors with 2.83 GHz and 8 Gbyte of memory (using only one
processor).

\subsection{Approximation by nonnegative splines}\label{ex:spl}
Consider the problem of approximating a one-dimensional function given
only by a large amount of noisy measurement by a cubic spline.
Additionally, we require that the function is nonnegative. This kind of
problem arises in many application, for instance, in shape optimization
considering unilateral contact or in arrival rate approximation
\cite{alizadeh}.

Assume that function $f:\RR\to \RR$ is defined on interval $[0,1]$. We
are given its function values $b_i$, $i=1,\ldots,n$ at points
$t_i\in(0,1)$. We may further assume that the function values are
subject to a random noise. We want to find a smooth approximation of
$f$ by a cubic spline, i.e., by a function of the form
\begin{equation}\label{eq:spl}
  P(t) = P^{(i)}(t) = \sum_{k=0}^3 P^{(i)}_k(t-a_{i-1})^k
\end{equation}
for a point $t\in[a_{i-1},a_i]$, where $0=a_0<a_1<\ldots < a_m=1$ are
the knots and $P^{(i)}_k (i=1,\ldots,m,\ k=0,1,2,3)$ the coefficients
of the spline. The spline property that $P$ should be continuous and
have continuous first and second derivatives is expressed by the
following equalities for $i=1,\ldots,m-1$:
\begin{align}
P^{(i+1)}_0 &\!- P^{(i)}_0 \!- P^{(i)}_1(a_i-a_{i-1}) - P^{(i)}_2(a_i-a_{i-1})^2 -
 P^{(i)}_3(a_i-a_{i-1})^3  =0\label{eq:spl1a}\\
P^{(i+1)}_1 &- P^{(i)}_1 - 2P^{(i)}_2(a_i-a_{i-1}) - 3P^{(i)}_3(a_i-a_{i-1})^2 =0\label{eq:spl1b}\\
2P^{(i+1)}_2 &- 2P^{(i)}_2 - 6P^{(i)}_3(a_i-a_{i-1})  =0\,. \label{eq:spl1c}
\end{align}
The function $f$ will be approximated by $P$ in the least square sense,
so we want to minimize
$$
  \sum_{j=1}^n (P(t_j) - b_j)^2
$$
subject to (\ref{eq:spl1a}),(\ref{eq:spl1b}),(\ref{eq:spl1c}).

Now, the original function $f$ is assumed to be nonnegative and we also
want the approximation $P$ to have this property. A simple way to
guarantee nonnegativity of a spline is to express is using $B$-splines
and consider only nonnegative $B$-spline coefficients. However, it was
shown by de Boor and Daniel \cite{deboor} that this may lead to a poor
approximation of $f$. In particular, they showed that while
approximation of a nonnegative function by nonnegative splines of order
$k$ gives errors of order $h^k$, approximation by a subclass of
nonnegative splines of order $k$ consisting of all those whose
$B$-spline coefficients are nonnegative may yield only errors of order
$h^2$. In order to get the best possible approximation, we use a result
by Nesterov \cite{nesterov00} saying that $P^{(i)}(t)$ from
(\ref{eq:spl}) is nonnegative if and only if there exist two symmetric
matrices
$$
  X^{(i)}=\begin{pmatrix}x_i&y_i\\y_i&z_i\end{pmatrix},\qquad
  S^{(i)}=\begin{pmatrix}s_i&v_i\\v_i&w_i\end{pmatrix}
$$
such that
\begin{align}
P^{(i)}_0 &= (a_i-a_{i-1})s_i\label{eq:spl2a}\\
P^{(i)}_1 &= x_i - s_i + 2(a_i-a_{i-1})v_i\label{eq:spl2b}\\
P^{(i)}_2 &= 2y_i - 2v_i + (a_i-a_{i-1})w_i \label{eq:spl2c}\\
P^{(i)}_3 & = z_i - w_i \label{eq:spl2d}\\
X^{(i)} & \succeq 0,\quad  S^{(i)} \succeq 0 \,.\label{eq:spl2e}
\end{align}

Summarizing, we want to solve an NLP-SDP problem
\begin{align}
&\min_{\stackrel{\scriptstyle P^{(i)}_k\in\RR}{\scriptstyle i=1,\ldots,m,\; k=0,1,2,3}}
\sum_{j=1}^n (P(t_j) - b_j)^2\label{spli}\\
&\mbox{subject to}\nonumber\\
&\qquad  (\ref{eq:spl1a}),(\ref{eq:spl1b}),(\ref{eq:spl1c}),\quad i=1,\ldots,m\nonumber\\
&\qquad (\ref{eq:spl2a})-(\ref{eq:spl2e}),\quad i=1,\ldots,m\,.\nonumber
\end{align}

More complicated (``more nonlinear'') objective functions can be
obtained when considering, for instance, the problem of approximating
the arrival rate function of a non-homogeneous Poisson process based on
observed arrival data \cite{alizadeh}.

\example A problem of approximating a cosine function given at 500
points by noisy data of the form
\verb!cos(4*pi*rand(500,1))+1+.5.*rand(500,1)-.25!
approximated by a nonnegative cubic spline with 7 knots lead to an NSDP
problem in 80 variables, 16 matrix variables, 16 matrix constraints,
and 49 linear inequality constraints. The problem was solved by \pennon
in about 1 second using 17 global and 93 Newton iterations.

\subsection*{Acknowledgements}
The authors would like to thank Didier Henrion and Johan L\"ofbeg for
their constant help during the code development. The work has been
partly supported by grant A100750802 of the Czech Academy of Sciences
(MK) and by DFG cluster of excellence 315 (MS). The manuscript was
finished while the first author was visiting the Institute for Pure and
Applied Mathematics, UCLA. The support and friendly atmosphere of the
Institute are acknowledged with gratitude.

\bibliographystyle{plain}
\bibliography{pennon}
%


\printindex

\end{document}